\documentclass[12pt,a4paper]{article}
\usepackage[utf8]{inputenc}
\usepackage{amsmath}
\usepackage{amsfonts}
\usepackage{amssymb}
\usepackage{makeidx}
\usepackage{graphicx}

\usepackage[left=2.25cm, right=2.25cm, top=2.25cm, bottom=2.25cm]{geometry}

\usepackage{moreverb}
\usepackage{bbm}

\newtheorem{d1}{Definition}
\newtheorem{d2}[d1]{Definition}

\newtheorem{p1}[d1]{Property}

\DeclareMathOperator\erfc{erfc}

\providecommand{\keywords}[1]
{
	\small	
	\textbf{\textit{Keywords---}} #1
}

\title{Abstract, keywords and references template}
\author{M. B\l asik$^{1}$  \\
	\small $^{1}$Institute of Mathematics, Czestochowa University of Technology \\
	Armii Krajowej 21, 42-201 Czestochowa, Poland \\
	marek.blasik@gmail.com
}
\date{} 

\title{A numerical method for the solution of the two-phase fractional Lam\'{e}-Clapeyron-Stefan problem}

\begin{document}
\maketitle

\abstract{In this paper we present a numerical solution of a two-phase fractional Stefan problem with time 
	derivative described in the Caputo sense. In the proposed algorithm, we use a special case of front-fixing method
	supplemented by the iterative procedure, which allows us to determine the position of the moving boundary. In 
	the final part, we also present some examples illustrating the comparison of the analytical solution with the 
	results received by  the proposed numerical method.}

\keywords{moving boundary problems, fractional derivatives and integrals, Stefan problems,  phase changes,  
	finite difference methods}

\section{Introduction}
Moving boundary problems ares one of the most important area within partial differential equations. This particular kind of boundary-value problem was originally intended to describe the solid-liquid phase change process, but also refers to such phenomena as solute transport, molecular diffusion or controlled drug release \cite{Hig61,Hig63,Coh88}. Moving boundary problems, as the name implies, are characterized by having a moving boundary of the domain, which has to be determined as a part of the solution. These kinds of problems are also called Stefan problems, in connection with the early work of Slovene scientist Joseph Stefan, who investigated ice formation in the polar Arctic seas \cite{Ste91}. In recent years, the classical Stefan problem has been very well studied and described in many papers and monographs (compare \cite{Cra84,Gup03,Hil87,Rub71,Ozi93} and the references therein). 

The fractional Stefan problem is a natural generalization of th classical Stefan problem. 
The first paper devoted to this issue was published in 2004 and concerned the mathematical modeling of the controlled release of a drug from slab matrices \cite{Liu04}.
Recently, there has been an increase in the number of scientific publications concerning the moving 
boundary problems modeled by the anomalous diffusion equation. Basically,
published papers deal with three classes of phenomena. As mentioned earlier, the first
concerned the controlled release of a drug from slab matrices \cite{Li09,Chen11}.
A second class of problems relates to mathematical modeling of the thermal conductivity
with phase transitions \cite{Sin11,Vol10,Vol13,Vol14} and the last class refers to mathematical modeling a movement of the shoreline in a sedimentary ocean basin 
\cite{Raj13,Raj13b}.

Most of the analytical solutions of fractional Stefan problem were obtained for the one-phase case.
Liu Junyi and Xu Mingyu \cite{Liu04} studied the one-phase Stefan problem with fractional
anomalous diffusion (Riemann-Liouville derivative with respect to time variable was used) and got 
an exact solution (concentration of the drug in the matrix) in terms of the Wright's function. They 
also showed that position of the penetration of solvent at time $t$ moving
as $\sim t^{\frac{\alpha}{2}}$, $0<\alpha\leq 1$. Xicheng Li et al. \cite{Xic08} considered
one-phase fractional Stefan problem with Caputo derivative with respect to time and two 
types of space-fractional derivatives (Caputo and Riemann-Liouville). They got the solution 
in terms of the generalized Wright's function. It should be noted that in both cases
deliberated by the authors, function describing the moving boundary is 
$\sim t^\frac{\alpha}{\beta}$, $0<\alpha\leq 1$, $1<\beta\leq 2$, where $\alpha$, 
$\beta$ denotes the orders of the fractional derivatives with respect to time and spatial 
variable respectively. In the paper \cite{Vol10} Voller analytically solved a limit case fractional 
Stefan problem describing the melting process. For the governing equation with a Caputo derivative with 
respect to time of order $0<\beta\leq 1$ and for the same fractional derivative with respect to space 
for the flux of order $0<\alpha\leq 1$, he showed that a melting front is described by power function 
$t^{\frac{\beta}{\alpha+1}}$. More analytical results were published in the papers \cite{Liu09,Ros13,Ros16}.

An important result with respect to this paper is the closed analytical solution  of the 
two-phase fractional Lam\'{e}-Clapeyron-Stefan problem in a semi-infinite region obtained by 
Roscani and Tarzia \cite{Ros14}, which will be recalled in Section 3. The problem involves 
determination of three functions, namely, $u_1$, $u_2$ and $S$ fulfilling two subdiffusion 
equations (\ref{eq:sub1n}), (\ref{eq:sub2n}) and additional differential equation (interface 
energy balance condition) (\ref{eq:ste}) governing function $S$. The authors showed that 
$u_1$ and $u_2$ can be expressed in terms of the Wright's function, moreover location 
of the phase-change interface is described by power function (\ref{eq:s}), where parameter $p$
can be evaluated form transcendental equation (\ref{eq:trans}). This solution will be especially 
useful for validation of the numerical method proposed in Section 4.

An alternative to the closed analytical solutions are those received by numerical methods. 
There are many techniques of numerical solving of classical Stefan problem, some of them have 
been generalized to the case of the fractional order. Xiaolong Gao et al. \cite{Xia15} 
generalized the boundary immobilisation technique (also known as front-fixing method  
\cite{Cra84}) to the fractional Stefan problem with a space-fractional derivative. Another variant 
of the front-fixing method was proposed in paper \cite{Bla15} and applied to the fractional 
Stefan problem with time-fractional derivative.

Our aim is to develop the numerical method of solving the two-phase, one-dimensional fractional 
Stefan problem. The proposed numerical scheme is an extension of the front-fixing method 
\cite{Bla15} to the two-phase problem. Our new approach is based on the suitable selection of 
the new space coordinates. The original coordinate system $(x,\tau)$ is transformed into a two 
new orthogonal systems $(v_1,\tau)$ and $(v_2,\tau)$ using transformation (\ref{eq:v1}) and
(\ref{eq:v2}), respectively. Both new spatial variables $v_1$, $v_2$ depend on the parameter $p$
which is unknown and chosen a priori. The proposed numerical method uses integro-differential 
equations equivalent to the corresponding governing differential equations of the problem.
The solutions of the integro-differential equations are obtained separately for each phase and 
fulfill the interface energy balance condition (the moving boundary is fixed ) only when the 
value of parameter $p$ is correct. Selection of the appropriate value of parameter $p$ is 
implemented by iterative algorithm on the basis of the fractional Stefan condition.

The paper is organized as follows. In the next Section, we introduce definitions of the fractional
integrals and derivatives together with some of their properties. In Section 3, we formulate
a mathematical model describing the melting process and recall closed analytical
similarity solution for two-phase fractional Stefan problem. Section 4 is devoted to the new 
numerical method of solving of two-phase fractional Stefan problem, which is a extension 
of the method developed in the paper \cite{Bla15} to the two-phase case. Section 5 contains the 
analytical and numerical results. The last Section includes a summary of the paper and 
conclusions.

\section{Preliminaries}
Let us recall the basic definitions and properties from fractional calculus
\cite{Kil06,Sam93} which will be further applied to formulate and solve the two-phase fractional
Stefan problem modeled by two linear equations with Caputo time derivative of order $\alpha \in
(0,1]$. First, we define the left-sided Riemann-Liouville integral.
\begin{d1}
	The left-sided Riemann-Liouville integral of order $\alpha$, denoted as $I_{0+}^{\alpha}$, is given by the following formula for $Re(\alpha)>0$:
	\begin{equation}
	I_{0+}^{\alpha}f(t):=\frac{1}{\Gamma(\alpha)}\int_{0}^{t}\frac{f(u)du}
	{(t-u)^{1-\alpha}},
	\end{equation}
	where $\Gamma$ is the Euler gamma function.
\end{d1}
The left-sided Caputo derivative of order $\alpha \in (0,1)$ denoted by ${}^{c}D_{0+}^{\alpha}$ is defined via the above left-sided Riemann-Liouville integral
${}^{c}D_{0+}^{\alpha} f(t):=I_{0+}^{1-\alpha}f'(t)$.
\begin{d2}
	Let $Re(\alpha)\in(0,1]$. The left-sided Caputo derivative of order $\alpha$
	is given by the formula:
	\begin{equation}
	{}^{c}D_{0+}^{\alpha}f(t):=\left\lbrace
	\begin{array}{ll}
	\frac{1}{\Gamma(1-\alpha)}
	\int_{0}^{t}\frac{f^{'}(u)du}{(t-u)^{\alpha}},&0<\alpha<1,\\
	\frac{d f(t)}{dt}, & \alpha=1.
	\end{array}
	\right.
	\label{eq:Cap}
	\end{equation}
\end{d2}
An important property of fractional integral operators given in Definition 1 
is called a semigroup property, which allows composition of two left-sided fractional 
integrals.
\begin{p1}[cf. Lemma~2.3 \cite{Kil06}]
	\label{pr:3}
	If $Re(\alpha)>0$, and $Re(\beta)>0$, then the equation
	\begin{equation}
	I_{0+}^{\alpha}I_{0+}^{\beta}f(t)=I_{0+}^{\alpha+\beta}f(t)
	\end{equation}
	is satisfied at almost every point $t\in[0,b]$ for $f(t)\in L_p(0,b)$ where 
	$1 \leq p\leq \infty$. If $\alpha+\beta>1$, then the above relation holds
	at any point of $[0,b]$.
\end{p1}
The composition rule of the left-sided Riemann-Liouville integral with the left-sided
Caputo derivative is a consequence of the semigroup property for fractional integrals.
\begin{p1}[cf. Lemma~2.22, \cite{Kil06}]
	\label{pr:4}
	Let function $f\in C^1(0,b)$. Then, the composition rule for the left-sided Riemann-Liouville integral
	and the left-sided Caputo derivative is given as follows:
	\begin{equation}
	I_{0+}^{\alpha}{}^c D_{0+}^{\alpha}f(t)=f(t)-f(0).
	\end{equation}
\end{p1}
The two-parameter Wright function determined in a complex plane plays important role in the theory of partial differential 
equations of fractional order. It is a generalization of the complementary error function.
\begin{d2}
	Let $\gamma > -1$, $\delta\in \mathbb{C}, z \in \mathbb{C}$. The  two-parameter Wright function is given as the 
	following series:
	\begin{equation}
	\label{wright}
	W(z;\gamma,\delta):=\sum_{k=0}^{\infty}\frac{z^k}{k!\Gamma(\gamma k+\delta)}.
	\end{equation}
\end{d2}
We note that for $\gamma=-\frac{1}{2}$, $\delta=1$ the above two-parameter Wright function
becomes complementary error function \cite{Sha15}:
\begin{equation}
W\left( -z; -\frac{1}{2},1\right)  =\erfc\left(\frac{z}{2} \right).
\end{equation}
For $\gamma=-\frac{1}{2}$, $\delta=-\frac{1}{2}$ two-parameter Wright function
can be expressed by the formula \cite{Pod99}:
\begin{equation}
W\left( -z; -\frac{1}{2},-\frac{1}{2}\right)  =\frac{1}{\sqrt{\pi}}\exp\left(-\frac{z^2}{4} \right).
\end{equation}

\section{Mathematical formulation of the problem}
\begin{figure}[h!t]
	\begin{center}
		\includegraphics[scale=1]{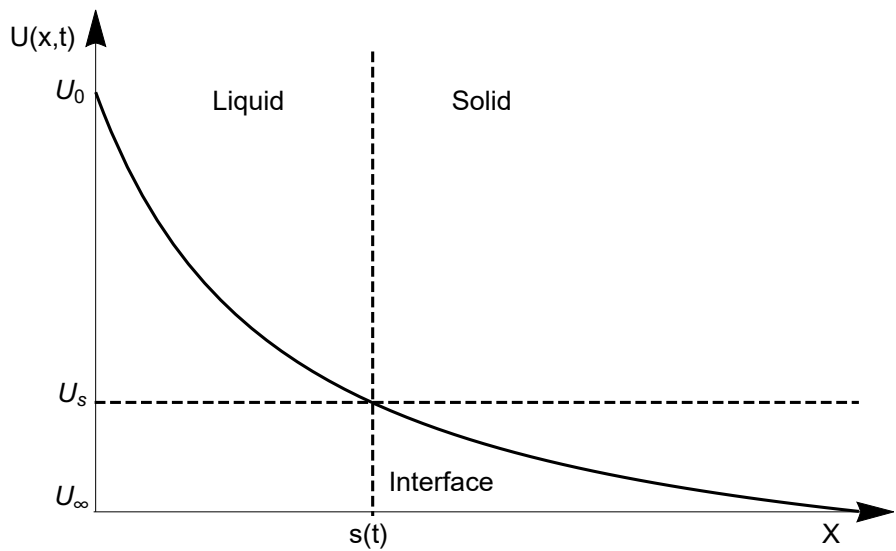}
		\caption{Semi-infinite slab melting from $x=0$ due to the high temperature $U_0$.}
		\label{fig:scheme}
	\end{center}
\end{figure}

The two-phase fractional Stefan problem is a mathematical model describing the
solidification and melting process. The mathematical formulation of the problem involves 
an anomalus diffusion equation for the liquid and solid phases and a condition at the liquid-solid
interface, called the fractional Stefan condition, which describes the position of the phase change
front. At the moving boundary $X=s(t)$, the temperature \textit{U} is constant and equal to 
melting point $U_s$. We are considering the melting of a semi-infinite, one dimensional slab occupying
$X\geq 0$, where the phase change interface moves from the heat source at a temperature
$U=U_0$ at the boundary $X=0$. A simple scheme of the model is shown in Figure
\ref{fig:scheme}.

Consider the following governing equations of the model
\begin{equation}
\label{eq:sub1}
\begin{split}
c_1\rho_1{}^{c}D_{0+,t}^{\alpha}U_1(X,t)=K_{1} \frac{\partial^2 U_1(X,t) }{\partial X^2},\quad 0<X<s(t)\leq l_1,\quad t>0,
\end{split}
\end{equation}

\begin{equation}
\label{eq:sub2}
\begin{split}
c_2\rho_2{}^{c}D_{0+,t}^{\alpha}U_2(X,t)=K_{2} \frac{\partial^2 U_2(X,t) }{\partial X^2},\quad s(t)<X<\infty,\quad t>0,
\end{split}
\end{equation}
where $U_1$ denotes the temperature in the liquid phase, $U_2$ the temperature in the
solid phase, $s(t)$ represents the position of the moving boundary, $K$ the modified
thermal conductivity (has SI unit [$J\cdot s^{-\alpha}\cdot m^{-1}\cdot K^{-1}$]), $\rho$ the density, $c$ the
specific heat and subscripts $1$ and $2$
indicate liquid and solid phase, respectively. Equations (\ref{eq:sub1},\ref{eq:sub2})
should be supplemented by Dirichlet conditions:
\begin{equation}
\begin{split}
U_1(0,t)=U_0,~ U_1(s(t),t)=U_2(s(t),t)=U_s,\lim_{X\to\infty} U_2(X,t)=U_{\infty},~ t>0.
\end{split}
\end{equation}
At $t=0$ the liquid phase does not exist, so we use two initial conditions:
\begin{equation}
U_2(X,0)=U_{\infty}, \quad s(0)=0.
\end{equation}
The position of the moving boundary $s(t)$ is determined by the fractional Stefan
condition, which expresses the heat balance in the melting layer:
\begin{equation}
\begin{split}
L \rho_1 {}^c{}D_{0+,t}^{\alpha}s(t)=\left. K_{2}\frac{\partial U_2(X,t)}{\partial X}\right|_{X=s(t)}\left.-K_{1}\frac{\partial U_1(X,t)}{\partial X}\right|_{X=s(t)},
\end{split}
\end{equation}
where $L$ is the latent heat. Let us note that the above mathematical model depends on eight
parameters. To simplify of the studied problem, we introduce the following dimensionless
variables
\[
\begin{split}
x=\frac{X}{l_1},\quad \tau=t\left(\frac{K_0}{c_0 \rho_1 l_1^2}\right)^\frac{1}{\alpha},u_1=
\frac{U_1-U_s}{U_0-U_s},\quad u_2=\frac{U_2-U_s}{U_0-U_s},\quad S=\frac{s}{l_1},
\end{split}
\]
which reduces the number of free parameters to five and makes it possible study their impact on the
solution.

We can rewrite governing equations (\ref{eq:sub1},\ref{eq:sub2}) in a non-dimensional form
\begin{equation}
\label{eq:sub1n}
\begin{split}
{}^{c}D_{0+,\tau}^{\alpha}u_1(x,\tau)=\kappa_{1} \frac{\partial^2 u_1(x,\tau) }{\partial x^2},\quad 0<x<S(\tau)\leq 1, \quad \tau>0,
\end{split}
\end{equation}

\begin{equation}
\label{eq:sub2n}
\begin{split}
{}^{c}D_{0+,\tau}^{\alpha}u_2(x,\tau)=\kappa_{2} \frac{\partial^2 u_2(x,\tau) }{\partial x^2},\quad S(\tau)<x<\infty, \quad \tau>0,
\end{split}
\end{equation}
supplemented with the boundary conditions
\begin{equation}
\label{eq:bound}
\begin{split}
u_1(0,\tau)=1,~ u_1(S(\tau),\tau)=u_2(S(\tau),\tau)=0,\quad \lim_{x\to\infty}u_2\left(x,
\tau\right)=\frac{U_{\infty}-U_s}{U_0-U_s},~\tau>0,
\end{split}
\end{equation}
initial conditions
\begin{equation}
\label{eq:ini1}
u_2(x,0)=\frac{U_{\infty}-U_s}{U_0-U_s},
\end{equation}
\begin{equation}
\label{eq:ini2}
S(0)=0,
\end{equation}
and the fractional Stefan condition
\begin{equation}
\label{eq:ste}
\begin{split}
{}^c{}D_{0+,\tau}^{\alpha}S(\tau)=\left. \lambda_{2}\frac{\partial u_2(x,\tau)}{\partial x}
\right|_{x=S(\tau)}\left.-\lambda_{1}\frac{\partial u_1(x,\tau)}{\partial x}\right|_{x=S(\tau)}
\end{split}
\end{equation}
where 
$\kappa_{1}=(K_{1}/K_0)(c_0/c_1)$, $\kappa_{2}=(K_{2}/K_0)(c_0/c_2)$, 
$\lambda_{1}=\frac{(U_0-U_s)K_{1}c_0}{L K_0}$, 
$\lambda_{2}=\frac{(U_0-U_s)K_{2}c_0}{L K_0}$,
are standard values of the respective variables.

According to results obtained by Roscani and Tarzia \cite{Ros14} the closed analytical solution of
the two-phase fractional Stefan problem (\ref{eq:sub1n})-(\ref{eq:ste}) is given by the functions
\begin{equation}
\label{eq:ph1}
u_1(x,\tau)=1-\frac{W\left (\frac{-x}{\sqrt{\kappa_{1}}\tau^{\alpha/2}};
	-\frac{\alpha}{2}, 1\right)-1}{W\left (\frac{-p}{\sqrt{\kappa_{1}}}; -\frac{\alpha}{2}, 1\right)-1},
\end{equation}
\begin{equation}
\label{eq:ph2}
\begin{split}
u_2(x,\tau)=\frac{U_{\infty}-U_s}{U_0-U_s} \frac{W\left (\frac{-p}
	{\sqrt{\kappa_{2}}}; -\frac{\alpha}{2}, 1\right)-W\left (\frac{-x}
	{\sqrt{\kappa_{2}}\tau^{\alpha/2}};-\frac{\alpha}{2}, 1\right)}
{W\left (\frac{-p}{\sqrt{\kappa_{2}}}; -\frac{\alpha}{2}, 1\right)},
\end{split}
\end{equation}
\begin{equation}
\label{eq:s}
S(\tau)=p \tau^{\alpha/2},
\end{equation}
\begin{equation}
\label{eq:trans}
\begin{split}
p\frac{\Gamma(1+\frac{\alpha}{2})}{\Gamma(1-\frac{\alpha}{2})}=
\frac{\lambda_{2}}{\sqrt{\kappa_{2}}}
\frac{U_{\infty}-U_s}{U_0-U_s}
\frac{W\left (\frac{-p}{\sqrt{\kappa_{2}}};-\frac{\alpha}{2}, 1-\frac{\alpha}{2}\right)}
{W\left (\frac{-p}{\sqrt{\kappa_{2}}}; -\frac{\alpha}{2}, 1\right)}
-\frac{\lambda_{1}}{\sqrt{\kappa_{1}}}\frac{W\left (\frac{-p}{\sqrt{\kappa_{1}}};-\frac{\alpha}{2}, 1-\frac{\alpha}{2}\right)}{W\left (\frac{-p}{\sqrt{\kappa_{1}}}; -\frac{\alpha}{2}, 1\right)-1}.
\end{split}
\end{equation}
It should be noted that the above solution reduces to the fractional
one-phase problem \cite{Liu09,Ros13} for $u_2\to 0$ .

When $\alpha\to 1$, then from (\ref{eq:ph1}-\ref{eq:trans}) the following
results of the classical two-phase Stefan problem are recovered:
\begin{equation}
u_1(x,\tau)=1-\frac{\erfc\left(\frac{x}{2\sqrt{\kappa_1 \tau}} \right) -1}
{\erfc\left(\frac{p}{2\sqrt{\kappa_1}} \right)-1 },
\end{equation}

\begin{equation}
u_2(x,\tau)=\frac{U_{\infty}-U_s}{U_0-U_s}\frac{\erfc\left(\frac{p}{2\sqrt{\kappa_2}}
	\right) -\erfc\left(\frac{x}{2\sqrt{\kappa_2 \tau}} \right) }
{\erfc\left(\frac{p}{2\sqrt{\kappa_2}}\right)},
\end{equation}

\begin{equation}
S(\tau)=p \sqrt{\tau},
\end{equation}

\begin{equation}
\frac{p}{2}=\lambda_2 \frac{U_{\infty}-U_s}{U_0-U_s}\frac{\exp\left(
	-\frac{p^2}{4\kappa_2}\right)}{\sqrt{\pi\kappa_2}\erfc\left(\frac{p}{2\sqrt{\kappa_2}}
	\right)}-\lambda_1 \frac{\exp\left(-\frac{p^2}{4\kappa_1}\right)}
{\sqrt{\pi\kappa_1}\left(\erfc\left(\frac{p}{2\sqrt{\kappa_1}} \right)-1 \right) }.
\end{equation}

\section{Numerical solution}
\begin{figure}[h!t]
	\begin{center}
		\includegraphics[scale=0.9]{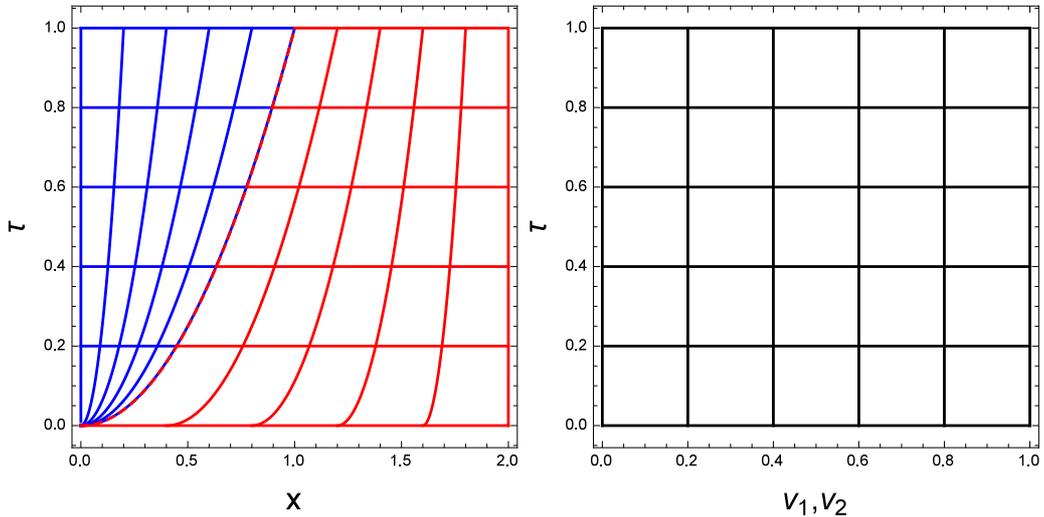}
		\caption{Meshes for $\alpha=1$, $p=1$, $m_1=m_2=5$, $n=5$, $\frac{l_2}{l_1}=2$.}
		\label{fig:mesh}
	\end{center}
\end{figure}

The numerical method proposed in the paper is an extension of the technique developed in \cite{Bla15} to the
two-phase problem. Applying the finite difference method to the fractional Stefan problem, we encounter a number of
difficulties related to the discretization of the fractional derivative with respect to the time variable. As we know, the
Caputo derivative is a non-local operator, which is defined on an interval. Suppose the positions, that the moving 
boundary will reach at different times are known and are at uniform spaced intervals in the same time layer, 
see at the left side of Figure \ref{fig:mesh}. For such a grid, the discretization of the Caputo derivative at some 
node with respect to time is very difficult because it requires values of function $u_1$ (or $u_2$) for all previous 
times (where the spatial variable is fixed) in points, which do not overlap with the mesh nodes. 
Let us also notice one important fact regarding discretization of equation (\ref{eq:sub1n}) leading 
to some ambiguity. The Caputo derivative 
requires for each point of the liquid phase their history from time $\tau = 0$, but they do not exist before they 
reach the melting point.  It seems that the only way to solve a mesh  problem is a suitable choice of new space 
coordinates for equations (\ref{eq:sub1n}) and (\ref{eq:sub2n}).

Let us first consider the region occupied by the first phase bounded by $x=0$ and $x=S(\tau)$, marked in blue in
Figure \ref{fig:mesh}. For the equation (\ref{eq:sub1n}), it is convenient to replace the spatial variable with the
following similarity variable \cite{Lan50}:
\begin{equation}
\label{eq:v1}
v_1=\frac{x}{S(\tau)}=\frac{x}{p\tau^{\alpha/2}},
\end{equation}
which has an important property, namely fixes the moving boundary at $v_1=1$ for all $\tau$.
We transform the first order derivative with respect to the spatial variable:
\begin{equation}
\frac{\partial u_1(x,\tau)}{\partial x}=\frac{\partial v_1}{\partial x}\frac{\partial u_1(v_1,\tau)}
{\partial v_1}=\frac{1}{p \tau^{\alpha/2}}\frac{\partial u_1(v_1,\tau)}{\partial v_1}.
\end{equation}
Subsequently, for the second order spatial derivative we have:
\begin{equation}
\label{eq:secder}
\frac{\partial^2 u_1(x,\tau)}{\partial x^2}=\frac{1}{p^2 \tau^{\alpha}}\frac{\partial^2 u_1(v_1,\tau)}{\partial v_1^2}.
\end{equation}
Respectively, for the first order partial derivative in time:
\begin{equation}
\frac{\partial u_1(x,\tau)}{\partial\tau}=\frac{\partial u_1(v_1,\tau)}{\partial\tau}-
\frac{\alpha v_1}{2\tau}\frac{\partial u_1(v_1,\tau)}{\partial v_1}.
\end{equation}
Caputo derivative of function $u_1(v_1,\tau)$ with respect to the time variable can be expressed as
follows:
\begin{equation}
\label{eq:caputo}
\begin{split}
&{}^c D_{0+,\tau}^{\alpha}u_1(v_1,\tau)=
\frac{1}{\Gamma(1-\alpha)}\int_{0}^{\tau}(\tau-\xi)^{-\alpha}
\frac{\partial}{\partial \xi} u_1(v_1,\xi)d\xi-\\&
\frac{\alpha v_1}{2\Gamma(1-\alpha)}\int_{0}^{\tau}(\tau-\xi)^{-\alpha}\frac{1}{\xi}
\frac{\partial}{\partial v_1} u_1(v_1,\xi)d\xi.
\end{split}
\end{equation}
Using formulas (\ref{eq:secder}) and (\ref{eq:caputo}) we write equation (\ref{eq:sub1n}) in 
coordinate system $(v_1,\tau)$:
\begin{equation}
\label{eq:sub}
\begin{split}
&\frac{1}{\Gamma(1-\alpha)}\int_{0}^{\tau}\frac{1}{(\tau-\xi)^{\alpha}}
\frac{\partial u_1(v_1,\xi)}{\partial \xi}d\xi-\frac{\alpha v_1}{2\Gamma(1-\alpha)}\int_{0}^{\tau}\frac{1}{(\tau-\xi)^{\alpha}}\frac{1}{\xi}
\frac{\partial u_1(v_1,\xi)}{\partial v_1}d\xi=\\&=
\frac{\kappa_{1}}{p^2\tau^{\alpha}}\frac{\partial^2 u_1(v_1,\tau)}{\partial v_1^2}
\end{split}
\end{equation}
We integrate the previous equation applying the left-sided Riemann-Liouville integral of order $\alpha\in(0,1)$ and we
get the following equation:
\begin{equation}
\label{eq:sub3}
\begin{split}
I_{0+,\tau}^{\alpha}{}^{c}D_{0+,\tau}^{\alpha}u_1(v_1,\tau)-\frac{\alpha 
	v_1}{2}I_{0+,\tau}^{\alpha} I_{0+,\tau}^{1-\alpha}\left(\frac{\partial}{\partial v_1}
\frac{u_1(v_1,\tau)}{\tau}\right)=\frac{\kappa_{1}}{p^2}I_{0+,\tau}^{\alpha}
\left(\frac{1}{\tau^{\alpha}}\frac{\partial^2 u_1(v_1,\tau)}{\partial v_1^2}\right)
\end{split}
\end{equation}
Finally, using Properties \ref{pr:3} and \ref{pr:4} to equation (\ref{eq:sub3}), we obtain an 
integro-differential equation in the form of:
\begin{equation}
\label{eq:fin1}
\begin{split}
u_1(v_1,\tau)&=u_1(v_1,0)+\frac{\alpha v_1}{2}\int_{0}^{\tau}
\frac{\partial}{\partial v_1}\frac{u_1(v_1,\xi)}{\xi}d\xi+
\frac{\kappa_{1}}{p^2 \Gamma(\alpha)}\int_{0}^{\tau}\frac{1}{(\tau-\xi)^{1-\alpha}}
\frac{1}{\xi^{\alpha}}\frac{\partial^2 u_1(v_1,\xi)}{\partial v_1^2}d\xi
\end{split}
\end{equation}
The kernel of the second integral on the right-hand side of formula (\ref{eq:fin1}) causes some 
difficulties in deriving a numerical scheme. For this reason we are introducing an auxiliary function:
\begin{equation}
\bar{u}_1(v_1,\tau)=u_1(v_1,\tau)\tau^{-\alpha},
\end{equation}
which leads to the integro-differential equation:
\begin{equation}
\label{eq:fin2}
\begin{split}
&\bar{u}_1(v_1,\tau)\tau^{\alpha}=\bar{u}_1(v_1,0)\tau_0^{\alpha}+\frac{\alpha v_1}{2}\int_{0}^{\tau}
\frac{\partial}{\partial v_1}\frac{\bar{u}_1(v_1,\xi)}{\xi^{1-\alpha}}d\xi+\\&+
\frac{\kappa_{1}}{p^2 \Gamma(\alpha)}\int_{0}^{\tau}\frac{1}{(\tau-\xi)^{1-\alpha}}
\frac{\partial^2 \bar{u}_1(v_1,\xi)}{\partial v_1^2}d\xi,
\end{split}
\end{equation}
supplemented with the boundary conditions
\[
\bar{u}_1(0,\tau)=\tau^{-\alpha}, \qquad \bar{u}_1(1,\tau)=0,
\]
and initial condition
\[
\bar{u}_1(v_1,0)=0.
\]

The region marked in blue in Figure \ref{fig:mesh} in the $(x,\tau)$ plane is transformed into the unit square
(in the general case to the rectangle) using the transformation (\ref{eq:v1}). Rectangular mesh created in this 
way consists of horizontal lines spaced $\Delta \tau=1/(np^{2/\alpha})$ units apart and vertical lines spaced 
$\Delta v_1=1/m_1$ units apart. The points  $((v_1)_i, \tau_j)=(i\Delta v_1,j\Delta \tau)$ where
$i=0,...,m_1$ and $j=0,...,n$, of intersection of the horizontal and vertical lines are grid nodes.
We denote the value of function $\bar{u}_1$ at  point $(i\Delta v_1,j\Delta \tau)$ by $(\bar{u}_1)_{i,j}$.
At this stage, the value of  parameter $p$ is not known and chosen a  priori. Construction
of similarity variable (\ref{eq:v1}) requires an additional assumption for variable $\tau$. Let $\tau_0$ be a very small positive number.

The method of discretization of the integro-differential equation (\ref{eq:fin2})
is described in detail in the paper \cite{Bla15} and leads to an implicit scheme:
\begin{equation}
\label{eq:implicit2}
\begin{split}
&(\bar{u}_1)_{i-1,k+1}(-r^1_{k+1,k+1}+q^1_{i,k+1})+(\bar{u}_1)_{i,k+1}
(\tau_{k+1}^{\alpha}
+2r^1_{k+1,k+1})+(\bar{u}_1)_{i+1,k+1}(-r^1_{k+1,k+1}-q^1_{i,k+1})=\\&=(\bar{u}_1)_{i,0}
\tau_0^{\alpha}
+\sum_{j=0}^{k}r^1_{j,k+1}((\bar{u}_1)_{i-1,j}-2(\bar{u}_1)_{i,j}
+(\bar{u}_1)_{i+1,j})+\sum_{j=1}^{k}q^1_{i,j} ((\bar{u}_1)_{i+1,j}-(\bar{u}_1)_{i-1,j}),
\end{split}
\end{equation}
where
\begin{equation}
r_{j,k+1}^1=\frac{c_{j,k+1} \kappa_{1}}{p^2 \Gamma(\alpha)(\Delta v_1)^2} 
,\qquad q_{i,j}^1=\frac{\alpha i \tau_j^{\alpha-1} \Delta\tau}{4},
\end{equation}
\begin{equation}
\label{eq:weight}
c_{j,k+1}=\left\{
\begin{array}{ll}
\frac{(\Delta\tau)^{\alpha}}{\alpha(\alpha+1)}\left( k^{\alpha+1}-(k-\alpha)(k+1)^{\alpha}
\right) & \textrm{for $j = 0$} \\
\frac{(\Delta\tau)^{\alpha}}{\alpha(\alpha+1)}( (k-j+2)^{\alpha+1}+(k-j)^{\alpha+1}&\\
-2(k-j+1)^{\alpha+1}) & \textrm{for $ 1\leq j \leq k$}\\
\frac{(\Delta\tau)^{\alpha}}{\alpha(\alpha+1)}& \textrm{for $j=k+1$}
\end{array} \right.
\end{equation}
Weights (\ref{eq:weight}) are the result of application of the trapezoidal rule \cite{Die10} to the last 
integral term in formula (\ref{eq:fin2}). Obtained numerical scheme can be written in the matrix form
\begin{equation}
\mathbf{A}\mathbf{(\bar{U}_1)}_{k+1}=\mathbf{B},
\end{equation}
where $\mathbf{\bar{U}_1}$ is a vector of unknown values of function $\bar{u_1}$ at  
instant $\tau_{k+1}$.
$$\mathbf{A}=\left[\begin{array}{cccccccc}
a_{k+1}^2 & a_{1,k+1}^3 & 0 & 0 & \cdots & 0 & 0 & 0 \\
a_{2,k+1}^1 & a_{k+1}^2 & a_{2,k+1}^3 & 0 & \cdots & 0 & 0 & 0 \\
0 & a_{3,k+1}^1 & a_{k+1}^2 & a_{3,k+1}^3 & \cdots & 0 & 0 & 0 \\
\vdots & \vdots & \vdots & \vdots & \ddots & \vdots & \vdots & \vdots\\
0 & 0 & 0 & a_{i,k+1}^1 & a_{k+1}^2 & a_{i,k+1}^3 & 0 & 0\\
\vdots & \vdots & \vdots & \vdots & \ddots & \vdots & \vdots & \vdots\\
0 & 0 & 0 & 0 & \cdots & a_{m-2,k+1}^1 & a_{k+1}^2 & a_{m-2,k+1}^3 \\
0 & 0 & 0 & 0 & \cdots & 0 & a_{m-1,k+1}^1 & a_{k+1}^2 \\
\end{array}
\right]_{(m_1-1)\times(m_1-1)}
$$
The elements of matrices $\mathbf{A}$ and $\mathbf{B}$ are defined as follows:
$$\mathbf{B}=\left[\begin{array}{c}
b_1-a_{1,k+1}^1 (\bar{u}_1)_{0,k+1}\\
b_2\\
b_3\\
\vdots\\
b_i\\
\vdots\\
b_{m_1-2}\\
b_{m_1-1}-a_{m_1-1,k+1}^3 (\bar{u}_1)_{m_1,k+1}
\end{array}
\right]_{(m_1-1)\times 1}
$$
where
\begin{eqnarray*}
	&&a_{i,k+1}^1:=-r^1_{k+1,k+1}+q^1_{i,k+1}\\
	&&a_{k+1}^2:=\tau_{k+1}^{\alpha}+2r^1_{k+1,k+1}\\
	&&a_{i,k+1}^3:=-r^1_{k+1,k+1}-q^1_{i,k+1}\\
	&&b_i:=(\bar{u}_1)_{i,0}\tau_0^{\alpha}+\sum_{j=0}^{k}r^1_{j,k+1}((\bar{u}_1)_{i-1,j}-2(\bar{u}_1)_{i,j}
	+(\bar{u}_1)_{i+1,j})
	+\sum_{j=1}^{k}q^1_{i,j} ((\bar{u}_1)_{i+1,j}-(\bar{u}_1)_{i-1,j}),
\end{eqnarray*}
We recover the values of function $u_1$ in coordinate system $(x,\tau)$ by applying
formulas
\begin{equation}
(u_1)_{i,j}=(\bar{u}_1)_{i,j}\tau_{j}^{\alpha},
\end{equation}
\begin{equation}
x_{i,j}=(v_1)_{i} p \tau_{j}^{\alpha/2}.
\end{equation}

The region marked in red in Figure \ref{fig:mesh} in the $(x,\tau)$ plane is transformed into the unit square using the following transformation:
\begin{equation}
\label{eq:v2}
v_2=\frac{x-p\tau^{\alpha/2}}{l_2/l_1-p\tau^{\alpha/2}},
\end{equation}
which  fixes the moving boundary at $v_2=0$ for all $\tau$.
Let us note, that the mesh for the semi-infinite region contains an infinite number of
nodes, which makes it impossible to perform any numerical calculations. Therefore, from a
practical point of view the Dirichlet boundary condition in infinity (\ref{eq:bound}) can
be replaced by the following condition:
\begin{equation}
u_2\left(l_2/l_1,\tau\right)=\frac{U_{\infty}-U_s}{U_0-U_s},
\end{equation}
where $l_2$ is a sufficiently large positive real number.

We transform the first order derivative with respect to the spatial variable:
\begin{equation}
\frac{\partial u_2(x,\tau)}{\partial x}=\frac{\partial v_2}{\partial x}\frac{\partial u_2(v_2,\tau)}
{\partial v_2}=\frac{1}{l_2/l_1-p \tau^{\alpha/2}}\frac{\partial u_2(v_2,\tau)}{\partial v_2}.
\end{equation}
Consequently, for the second-order spatial derivative
\begin{equation}
\label{eq:secder2}
\frac{\partial^2 u_2(x,\tau)}{\partial x^2}=\frac{1}{(l_2/l_1-p \tau^{\alpha/2})^2}\frac{\partial^2 u_2(v_2,\tau)}{\partial v_2^2}.
\end{equation}
Respectively, for the first order partial derivative with respect to the time variable:
\begin{equation}
\frac{\partial u_2(x,\tau)}{\partial\tau}=\frac{\partial u_2(v_2,\tau)}{\partial\tau}-
\frac{\alpha p (v_2-1)}{2\tau(p-l_2/l_1 \tau^{-\alpha/2})}\frac{\partial u_2(v_2,\tau)}{\partial v_2}.
\end{equation}
Caputo derivative of function $u_1(v_1,\tau)$ with respect to the time variable can be expressed as
follows:
\begin{equation}
\label{eq:caputo2}
\begin{split}
&{}^c D_{0+,\tau}^{\alpha}u_2(v_2,\tau)=
\frac{1}{\Gamma(1-\alpha)}\int_{0}^{\tau}(\tau-\xi)^{-\alpha}
\frac{\partial u_2(v_2,\xi)}{\partial \xi} d\xi-\\&
\frac{\alpha p (v_2-1)}{2\Gamma(1-\alpha)}\int_{0}^{\tau}
\frac{(\tau-\xi)^{-\alpha}}{\xi (p-l_2/l_1 \xi^{-\alpha/2})}
\frac{\partial u_2(v_2,\xi)}{\partial v_2} d\xi.
\end{split}
\end{equation}
Using formulas (\ref{eq:secder2}) and (\ref{eq:caputo2}) we can write equation (\ref{eq:sub2n}) in the new coordinate system $(v_2,\tau)$:
\begin{equation}
\label{eq:integro}
\begin{split}
&\frac{1}{\Gamma(1-\alpha)}\int_{0}^{\tau}\frac{1}{(\tau-\xi)^{\alpha}}
\frac{\partial u_2(v_2,\xi)}{\partial \xi}d\xi-
\frac{\alpha p (v_2-1)}{2\Gamma(1-\alpha)}\int_{0}^{\tau}
\frac{(\tau-\xi)^{-\alpha}}{\xi (p-l_2/l_1 \xi^{-\alpha/2})}
\frac{\partial u_2(v_2,\xi)}{\partial v_2} d\xi=\\&=
\frac{\kappa_{2}}{(l_2/l_1-p\tau^{\alpha/2})^2}\frac{\partial^2 u_2(v_2,\tau)}{\partial v_2^2}
\end{split}
\end{equation}
We integrate both sites of equation (\ref{eq:integro}) applying the left-sided
Riemann-Liouville integral of order $\alpha\in(0,1)$:
\begin{equation}
\label{eq:integro2}
\begin{split}
&I_{0+,\tau}^{\alpha}{}^{c}D_{0+,\tau}^{\alpha}u_2(v_2,\tau)-\frac{\alpha p 
	(v_2-1)}{2}I_{0+,\tau}^{\alpha} I_{0+,\tau}^{1-\alpha}\left(\frac{\partial}{\partial v_2}
\frac{u_2(v_2,\tau)}{\tau (p-l_2/l_1 \tau^{-\alpha/2})}\right)=\\&=I_{0+,\tau}^{\alpha}
\left(\frac{\kappa_{2}}{(l_2/l_1-p\tau^{\alpha/2})^2}\frac{\partial^2 u_2(v_2,\tau)}
{\partial v_2^2}\right).
\end{split}
\end{equation}
Finally, using Properties \ref{pr:3} and \ref{pr:4} to equation (\ref{eq:integro2}), we obtain an 
integro-differential equation in the form of:
\begin{equation}
\label{eq:fin12}
\begin{split}
u_2(v_2,\tau)&=u_2(v_2,0)+\frac{\alpha p (v_2-1)}{2}\int_{0}^{\tau}
\frac{\partial}{\partial v_2}\frac{u_2(v_2,\xi)}{\xi (p-l_2/l_1 \xi^{-\alpha/2})}d\xi+\\&+
\frac{\kappa_{2}}{\Gamma(\alpha)}\int_{0}^{\tau}\frac{(\tau-\xi)^{\alpha-1}}
{(l_2/l_1-p\xi^{\alpha/2})^2}
\frac{\partial^2 u_2(v_2,\xi)}{\partial v_2^2}d\xi.
\end{split}
\end{equation}
We introduce the second auxiliary function
\begin{equation}
\bar{u}_2(v_2,\tau)=\frac{u_2(v_2,\tau)}{(l_2/l_1-p\tau^{\alpha/2})^2},
\end{equation}
which leads to the integro-differential equation:
\begin{equation}
\label{eq:integro3}
\begin{split}
&\bar{u}_2(v_2,\tau)(l_2/l_1-p\tau^{\alpha/2})^2=
\bar{u}_2(v_1,0)(l_2/l_1-p\tau_0^{\alpha/2})^2
+\\&+\frac{\alpha p (v_2-1)}{2}\int_{0}^{\tau}
\frac{\partial}{\partial v_2}\left(p\xi^{\alpha-1}-\frac{l_2}{l_1}\xi^{\alpha/2-1}\right)
\bar{u}_2(v_2,\xi)d\xi+\\&+
\frac{\kappa_{2}}{\Gamma(\alpha)}\int_{0}^{\tau}\frac{1}{(\tau-\xi)^{1-\alpha}}
\frac{\partial^2 \bar{u}_2(v_2,\xi)}{\partial v_2^2}d\xi,
\end{split}
\end{equation}
supplemented with the boundary conditions

\begin{equation}
\bar{u}_2(0,\tau)=0, \qquad \bar{u}_2(1,\tau)=\frac{1}
{(l_2/l_1-p\tau^{\alpha/2})^2}\frac{U_{\infty}-U_s}{U_0-U_s},
\end{equation}
and initial condition
\begin{equation}
\bar{u}_2(v_2,0)=\frac{1}{(l_2/l_1-p\tau_0^{\alpha/2})^2}\frac{U_{\infty}-U_s}{U_0-U_s}.
\end{equation}

The integro-differential equation (\ref{eq:integro3}) is discretized in analogy to
equation (\ref{eq:fin2}). For the second phase we also operate on
a rectangular uniform grid consist of horizontal lines spaced $\Delta
\tau=1/(np^{2/\alpha})$ units apart and vertical lines spaced 
$\Delta v_2=1/m_2$ units apart. The points  $((v_2)_i, \tau_j)=(i\Delta v_2,j\Delta \tau)$ where $i=0,...,m_2$ and $j=0,...,n$, of intersection of the horizontal
and vertical lines are grid nodes. At this stage of the calculation, the value (for both
phases the same) of parameter $p$ is not known and chosen a  priori.

The implicit numerical scheme for second phase is in the form of:
\begin{equation}
\begin{split}
&(\bar{u}_2)_{i-1,k+1}(-r^2_{k+1,k+1}+q^2_{i,k+1})+(\bar{u}_2)_{i,k+1}
((l_2/l_1-p\tau_{k+1}^{\alpha/2})^2\\&
+2r^2_{k+1,k+1})+(\bar{u}_2)_{i+1,k+1}(-r^2_{k+1,k+1}-q^2_{i,k+1})=\\&=(\bar{u}_1)_{i,0}
(l_2/l_1-p\tau_0^{\alpha/2})^2
+\sum_{j=0}^{k}r^2_{j,k+1}((\bar{u}_2)_{i-1,j}-2(\bar{u}_2)_{i,j}
+\\&+(\bar{u}_2)_{i+1,j})+\sum_{j=1}^{k}q^2_{i,j} ((\bar{u}_2)_{i+1,j}-(\bar{u}_2)_{i-1,j}),
\end{split}
\end{equation}
where
\begin{equation}
r_{j,k+1}^2=\frac{c_{j,k+1} \kappa_{2}}{\Gamma(\alpha)(\Delta v_2)^2},\qquad
q_{i,j}^2=\frac{\alpha p (i \Delta v_2-1)(p\tau_j^{\alpha-1}-l_2/l_1 \tau_j^{\alpha/2-1})\Delta \tau}
{4 \Delta v_2}
\end{equation}
The obtained numerical scheme can be written in the matrix form
\begin{equation}
\mathbf{D}\mathbf{(\bar{U}_2)}_{k+1}=\mathbf{E},
\end{equation}
where $\mathbf{\bar{U}_2}$ is a vector of unknown values of function $\bar{u}_2$ at  
instant $\tau_{k+1}$.
$$\mathbf{D}=\left[\begin{array}{cccccccc}
d_{k+1}^2 & d_{1,k+1}^3 & 0 & 0 & \cdots & 0 & 0 & 0 \\
d_{2,k+1}^1 & d_{k+1}^2 & d_{2,k+1}^3 & 0 & \cdots & 0 & 0 & 0 \\
0 & d_{3,k+1}^1 & d_{k+1}^2 & d_{3,k+1}^3 & \cdots & 0 & 0 & 0 \\
\vdots & \vdots & \vdots & \vdots & \ddots & \vdots & \vdots & \vdots\\
0 & 0 & 0 & d_{i,k+1}^1 & d_{k+1}^2 & d_{i,k+1}^3 & 0 & 0\\
\vdots & \vdots & \vdots & \vdots & \ddots & \vdots & \vdots & \vdots\\
0 & 0 & 0 & 0 & \cdots & d_{m-2,k+1}^1 & d_{k+1}^2 & d_{m-2,k+1}^3 \\
0 & 0 & 0 & 0 & \cdots & 0 & d_{m-1,k+1}^1 & d_{k+1}^2 \\
\end{array}
\right]_{(m_2-1)\times(m_2-1)}
$$

$$\mathbf{E}=\left[\begin{array}{c}
e_1-d_{1,k+1}^1 (\bar{u}_2)_{0,k+1}\\
e_2\\
e_3\\
\vdots\\
e_i\\
\vdots\\
e_{m_2-2}\\
e_{m_2-1}-d_{m_2-1,k+1}^3 (\bar{u}_2)_{m_2,k+1}
\end{array}
\right]_{(m_2-1)\times 1}
$$
The elements of matrices $\mathbf{D}$ and $\mathbf{E}$ are defined as follows:
\begin{eqnarray*}
	&&d_{i,k+1}^1:=-r^2_{k+1,k+1}+q^2_{i,k+1}\\
	&&d_{k+1}^2:=(l_2/l_1-p \tau_{k+1}^{\alpha/2})^2+2r^2_{k+1,k+1}\\
	&&d_{i,k+1}^3:=-r^2_{k+1,k+1}-q^2_{i,k+1}\\
	&&e_i:=(\bar{u}_2)_{i,0}(l_2/l_1-p \tau_{0}^{\alpha/2})^2
	+\sum_{j=0}^{k}r^2_{j,k+1}((\bar{u}_2)_{i-1,j}-2(\bar{u}_2)_{i,j}
	+(\bar{u}_2)_{i+1,j})\\
	&&+\sum_{j=1}^{k}q^2_{i,j} ((\bar{u}_2)_{i+1,j}-(\bar{u}_2)_{i-1,j}),
\end{eqnarray*}
Applying the following formulas:
\begin{equation}
(u_2)_{i,j}=(\bar{u}_2)_{i,j}(l_2/l_1-p\tau_j^{\alpha/2})^2,
\end{equation}
\begin{equation}
x_{i,j}=(v_2)_{i}(l_2/ l_1 -p\tau_j^{\alpha/2})+p\tau_j^{\alpha/2},
\end{equation}
we recover the values of function $u_2$ in coordinate system $(x,\tau)$.

The presented numerical scheme allows us solve the governing equations 
(\ref{eq:sub1n}, \ref{eq:sub2n}), but only when the correct value of the parameter $p$ is
known. We can determine it using the fractional Stefan condition. First, we integrate both
sides of equation (\ref{eq:ste}) applying the left-sided Riemann-Liouville integral of
order $\alpha\in(0,1)$. Next, using Property \ref{pr:4} we obtain:
\begin{equation}
\label{eq:steint}
\begin{split}
S(\tau)-S(0)=\left. \lambda_{2}I_{0+,\tau}^{\alpha}\frac{\partial u_2(x,\tau)}{\partial x}
\right|_{x=S(\tau)}\left.-\lambda_{1}I_{0+,\tau}^{\alpha}
\frac{\partial u_1(x,\tau)}{\partial x}\right|_{x=S(\tau)}.
\end{split}
\end{equation}
Applying the difference quotient approximation for the first order derivative with respect to
space variable, weights (\ref{eq:weight}) and condition (\ref{eq:ini2}), we get:
\begin{equation}
\label{eq:Sfin}
S(\tau_n)=\frac{\lambda_{2}}{\Gamma(\alpha)}\sum_{j=0}^{n}c_{j,n}
\frac{(u_2)_{1,j}-(u_2)_{0,j}}{x_{1,j}-x_{0,j}}-
\frac{\lambda_{1}}{\Gamma(\alpha)}\sum_{j=0}^{n}c_{j,n}
\frac{(u_1)_{m_1,j}-(u_1)_{m_1-1,j}}{x_{m_1,j}-x_{m_1-1,j}},~j=0,...,n.
\end{equation}
Let us note that the value of function $S$ at the final time instant $\tau_n$ for the correct value of parameter $p$ should be equal to 1.
From this relation the  following criterion of convergence results in:
\begin{equation}
\label{eq:opt1}
|1-S(\tau_n,p)|<\epsilon,
\end{equation}
where $\epsilon>0$ is a certain arbitrarily small real number. We denote as $S(\tau_n,p)$ the value of $S(\tau_n)$ for a fixed value of parameter $p$. Below we present an iterative
algorithm for determining the value of parameter $p$, based on the bisection method
\cite{Bla15}.
\begin{enumerate}
	\item choose interval $[p_{a},p_{b}]$ for parameter $p$, define the value of parameter $\epsilon$
	\item calculate numerically  the solution of  equations (\ref{eq:fin2},
	\ref{eq:integro3}) for $p_{a}$ and $p_{b}$
	\item if one of the conditions:
	\\ $|1-S(\tau_n,p_{a})|<\epsilon$ or $|1-S(\tau_n,p_{b})|<\epsilon$
	\\ is fulfilled, then end the algorithm otherwise go to the  next step
	\item if the condition:\\
	$(1-S(\tau_n,p_{a}))\cdot(1-S(\tau_n,p_{b}))<0$\\
	is fulfilled, then go to the next step, otherwise go to step 1
	\item calculate $p_{c}:=\frac{p_{a}+p_{b}}{2}$
	\item calculate numerically  the solution of equations (\ref{eq:fin2},
	\ref{eq:integro3}) for $p_{c}$
	\item if $|1-S(\tau_n,p_{c})|<\epsilon$\\
	is fulfilled, then end the algorithm, otherwise go to the next step
	\item if $(1-S(\tau_n,p_{a}))\cdot(1-S(\tau_n,p_{c}))>0$\\
	then substitute $p_{a}:=p_{c}$ and go to the step 5 or\\
	if $(1-S(\tau_n,p_{b}))\cdot(1-S(\tau_n,p_{c}))>0$\\
	then substitute $p_{b}:=p_{c}$ and go to the step 5.
\end{enumerate}
The discussed algorithm was used to determine the value of the parameter $p$, whose values are presented in 
Table \ref{tab:2}.

\section{Numerical examples}
In order to validate the results obtained by the proposed numerical method with the
analytical solution one, twelve computer simulations were performed. We assumed four values
of order of the Caputo derivative $\alpha\in\{0.25,0.5,0.75,1\}$ and three sets of
physical parameter values. Due to the fact that the proposed method is iterative
(involves multiple solving of governing equations for different values of parameter $p$),
the following mesh parameters were adopted: $l_2/l_1=10$, $m_1=100$, $m_2=500$, $n=400$.

Table \ref{tab:1} shows the values of coefficient $p$ obtained from transcendental
equation (\ref{eq:trans}) depending on order of the Caputo derivative and three sets of physical parameters.
\begin{table}[h!t]
	\caption{Value of parameter \textit{p} from exact solution.}
	\begin{center}
		\begin{tabular}{|c|c|c|c|c|c|c|c|}
			\hline 
			\multicolumn{4}{|c|}{} &\multicolumn{4}{|c|}{$\alpha$}  \\ 
			\hline 
			$\lambda_{1}$ & $\lambda_{2}$ & $\kappa_{1}$ & $\kappa_{2}$
			& 0.25 & 0.5 & 0.75 & 1\\ 
			\hline \hline
			1 & 1 & 1 & 1 & 0.6834 & 0.7472 & 0.8299 & 0.9397 \\ 
			\hline
			1 & 2 & 1 & 1 & 0.5496 & 0.6013 & 0.668  & 0.7555 \\ 
			\hline 
			1 & 1 & 2 & 1 & 0.7218 & 0.7868 & 0.8697 & 0.9783 \\ 
			\hline 
		\end{tabular}
	\end{center}
	\label{tab:1}
\end{table}
The results collected in Table \ref{tab:1} indicates, that for a fixed set of physical
parameters, $p$ is an increasing function of order $\alpha$. This means that for
$\alpha$ values less than one the melting process is slower than in the case of the
classical Stefan problem. 

Let us now analyze the impact of the physical parameters on
the solution for the fixed value of $\alpha$. Comparing the first and second row of 
the Table \ref{tab:1}, we notice that $p$ is a decreasing function of $\lambda_{2}$.
This observation results directly from equation (\ref{eq:ste}). For a growing heat a flux
in solid zone ($\lambda_2$ is increasing), heat balance in the melting layer is preserved
only when value of the Caputo derivative of function $S$ decreases.
For $\alpha=1$ fractional derivative of function $S$ can be interpreted as a velocity of
movement of the melting front. The opposite behavior of the melting process is observed by
analyzing the first and third row of Table \ref{tab:1} for the fixed value of $\alpha$  
- coefficient $p$ is an increasing function of parameter $\kappa_{1}$.

Results given in Table \ref{tab:2} were obtained by applying the criterion of convergence
formulated in (\ref{eq:opt1}) and the iterative algorithm discussed in the previous
section. 

Table \ref{tab:3} contains the values of the variable $\tau$ for which the liquid phase 
propagates to the region bounded by $S(\tau)=1$. The collected data leads to the
following conclusion - changing parameter $\lambda_{2}$ has a greater effect on process
dynamics, than the changing parameter $\kappa_{1}$ for the fixed value of $\alpha$.
\begin{table}[h!t]
	\caption{Value of parameter \textit{p} from numerical solution.}
	\begin{center}
		\begin{tabular}{|c|c|c|c|c|c|c|c|}
			\hline 
			\multicolumn{4}{|c|}{} &\multicolumn{4}{|c|}{$\alpha$}  \\ 
			\hline 
			$\lambda_{1}$ & $\lambda_{2}$ & $\kappa_{1}$ & $\kappa_{2}$
			& 0.25 & 0.5 & 0.75 & 1\\ 
			\hline \hline
			1 & 1 & 1 & 1 & 0.7053 & 0.7358 & 0.8041 & 0.9311 \\ 
			\hline
			1 & 2 & 1 & 1 & 0.5691 & 0.5935 & 0.6497 & 0.7512 \\ 
			\hline 
			1 & 1 & 2 & 1 & 0.739  & 0.7801 & 0.8514 & 0.9674 \\ 
			\hline 
		\end{tabular}
	\end{center}
	\label{tab:2}
\end{table}
\begin{table}[h!t]
	\caption{Value of variable $\tau$ fulfilling equation $S(\tau)=1$ obtained from numerical solution.}
	\begin{center}
		\begin{tabular}{|c|c|c|c|c|c|c|c|}
			\hline 
			\multicolumn{4}{|c|}{} &\multicolumn{4}{|c|}{$\alpha$}  \\ 
			\hline 
			$\lambda_{1}$ & $\lambda_{2}$ & $\kappa_{1}$ & $\kappa_{2}$
			& 0.25 & 0.5 & 0.75 & 1\\ 
			\hline \hline
			1 & 1 & 1 & 1 & 16.329 & 3.411 & 1.789 & 1.153 \\ 
			\hline
			1 & 2 & 1 & 1 & 90.885 & 8.06 & 3.158 & 1.772 \\ 
			\hline 
			1 & 1 & 2 & 1 & 11.241  & 2.7 & 1.536 & 1.069 \\ 
			\hline 
		\end{tabular}
	\end{center}
	\label{tab:3}
\end{table}

\begin{figure}[h!t]
	\begin{center}
		\includegraphics[scale=1.4]{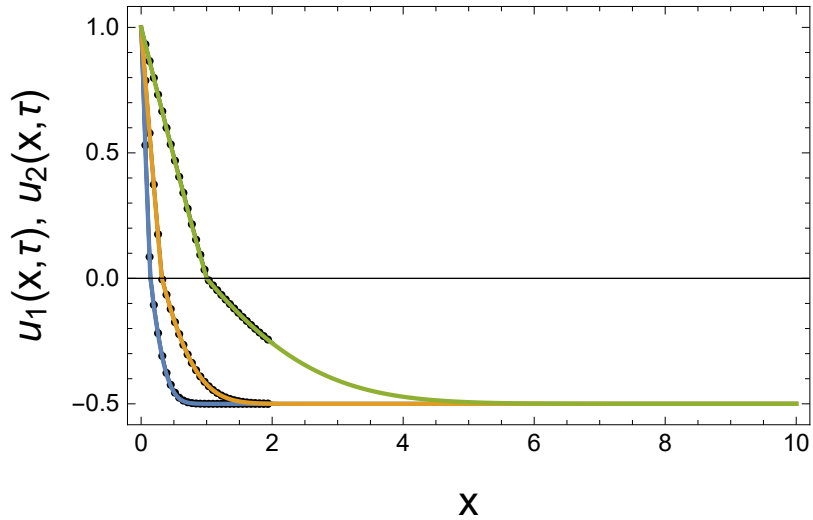}
		\caption{Numerical solution of $u_1$, $u_2$ for $\alpha=1$, $\lambda_1=1$, $\lambda_2=2$,
			$\kappa_{1}=1$, $\kappa_{2}=1$.}
		\label{fig:3}
	\end{center}
\end{figure}

\begin{figure}[h!t]
	\begin{center}
		\includegraphics[scale=1.4]{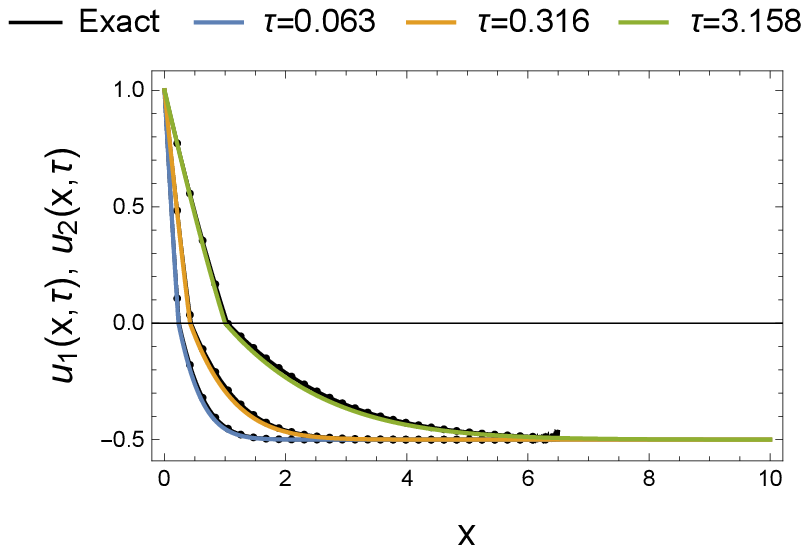}
		\caption{Numerical solution of $u_1$, $u_2$ for $\alpha=0.75$, $\lambda_1=1$, $\lambda_2=2$,
			$\kappa_{1}=1$, $\kappa_{2}=1$.}
		\label{fig:5}
	\end{center}
\end{figure}

\begin{figure}[h!t]
	\begin{center}
		\includegraphics[scale=1.4]{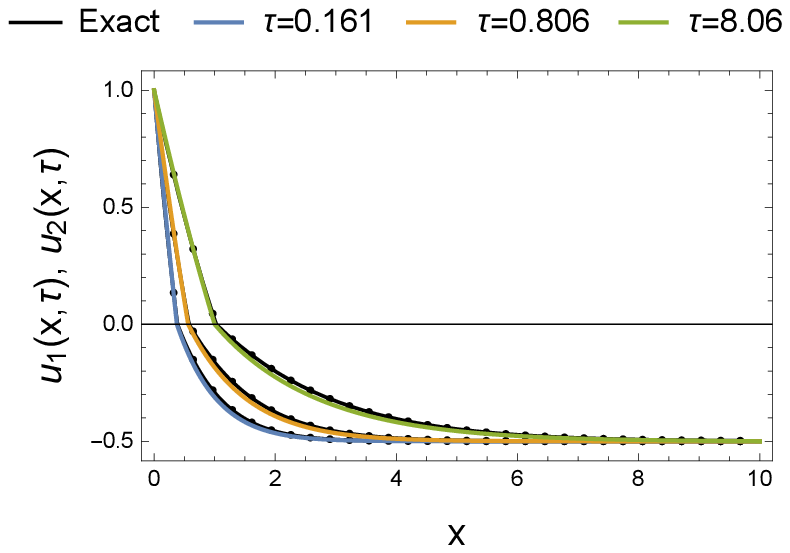}
		\caption{Numerical solution of $u_1$, $u_2$ for $\alpha=0.5$, $\lambda_1=1$, $\lambda_2=2$,
			$\kappa_{1}=1$, $\kappa_{2}=1$.}
		\label{fig:7}
	\end{center}
\end{figure}

\begin{figure}[h!t]
	\begin{center}
		\includegraphics[scale=1.4]{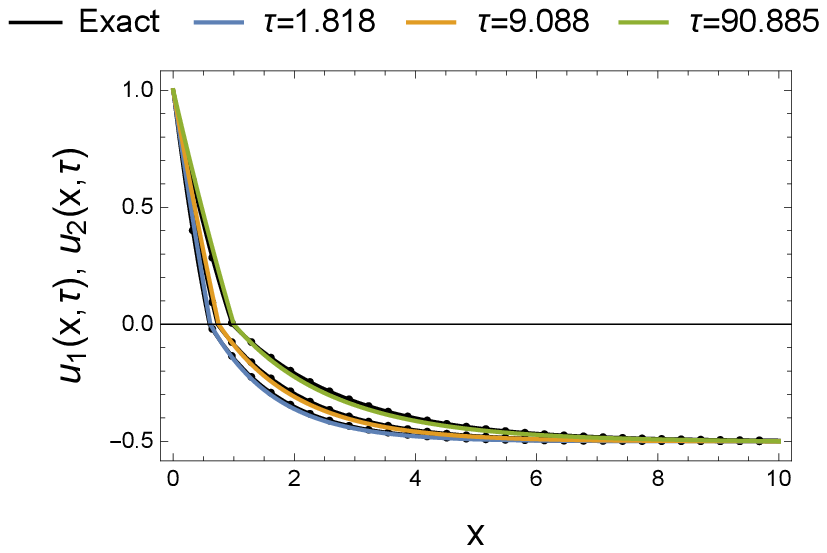}
		\caption{Numerical solution of $u_1$, $u_2$ for $\alpha=0.25$, $\lambda_1=1$, $\lambda_2=2$,
			$\kappa_{1}=1$, $\kappa_{2}=1$.}
		\label{fig:9}
	\end{center}
\end{figure}

In Figure \ref{fig:3},\ref{fig:5},\ref{fig:7} and \ref{fig:9} we present
dimensionless temperature profiles. On each graph we show the solutions obtained numerically 
(blue, yellow and green line) and the corresponding analytical solutions (black lines) plotted for
$x\in [0,10]$. In two cases for $\alpha=1$ and $\alpha=0.75$ we have adopted $x\in [0, 2]$ and 
$x\in [0, 6.5]$ respectively, because analytical solution (\ref{eq:ph2}) is difficult to calculate 
for large values of the spatial variable $x$. This is related to the definition of the Wright function, 
see formula (\ref{wright}). Therefore, due to time constraints of calculations, we used 500 terms 
of the series. However, this is not sufficient for large values of the similarity variable 
$\frac{x}{\tau^{\alpha/2}}$.

\begin{figure}[h!t]
	\begin{center}
		\includegraphics[scale=1.4]{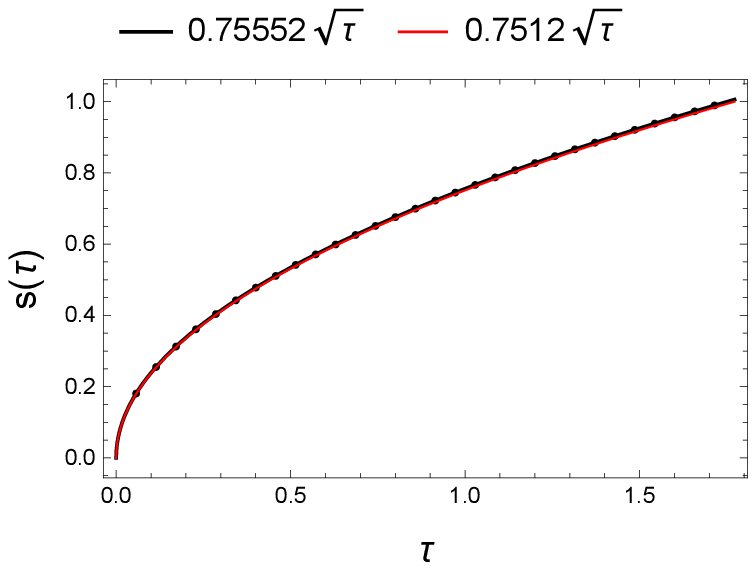}
		\caption{Numerical solution of $S$ for $\alpha=1$, $\lambda_1=1$, $\lambda_2=2$,
			$\kappa_{1}=1$, $\kappa_{2}=1$.}
		\label{fig:4}
	\end{center}
\end{figure}

\begin{figure}[h!t]
	\begin{center}
		\includegraphics[scale=1.4]{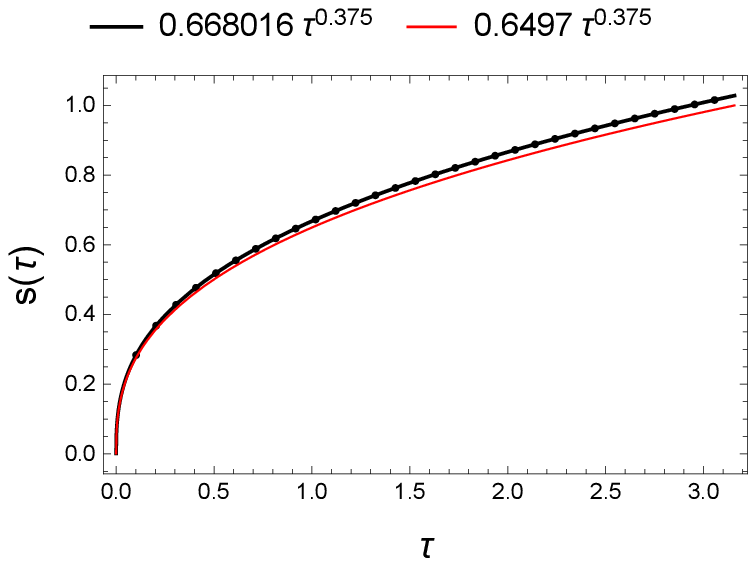}
		\caption{Numerical solution of $S$ for $\alpha=0.75$, $\lambda_1=1$, $\lambda_2=2$,
			$\kappa_{1}=1$, $\kappa_{2}=1$.}
		\label{fig:6}
	\end{center}
\end{figure}

\begin{figure}[h!t]
	\begin{center}
		\includegraphics[scale=1.4]{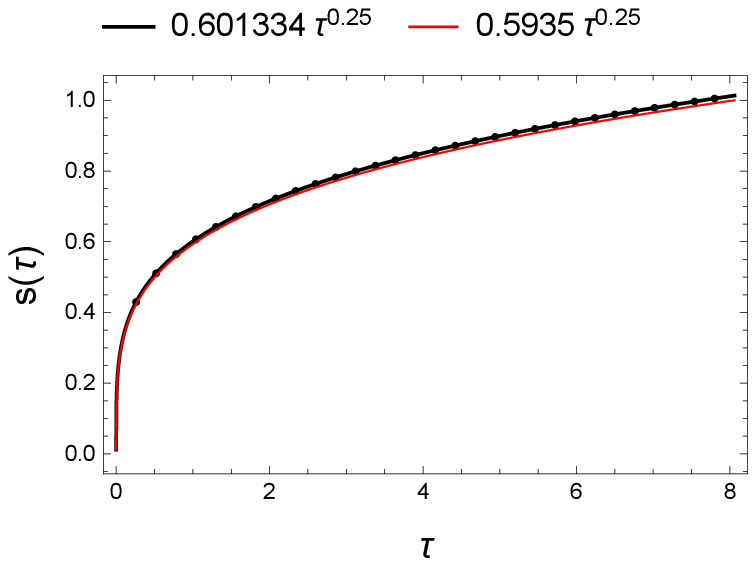}
		\caption{Numerical solution of $S$ for $\alpha=0.5$, $\lambda_1=1$, $\lambda_2=2$,
			$\kappa_{1}=1$, $\kappa_{2}=1$.}
		\label{fig:8}
	\end{center}
\end{figure}

\begin{figure}[h!t]
	\begin{center}
		\includegraphics[scale=1.4]{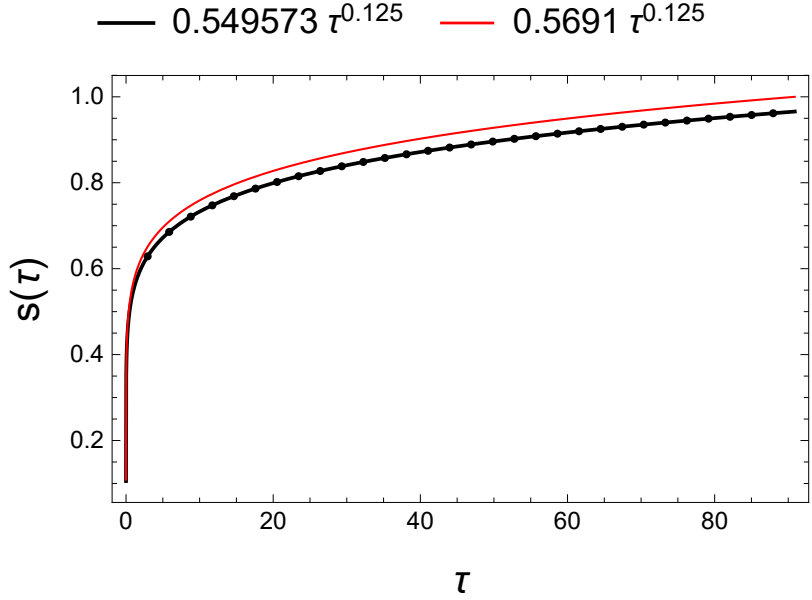}
		\caption{Numerical solution of $S$ for $\alpha=0.25$, $\lambda_1=1$, $\lambda_2=2$,
			$\kappa_{1}=1$, $\kappa_{2}=1$.}
		\label{fig:10}
	\end{center}
\end{figure}

In Figure \ref{fig:4},\ref{fig:6},\ref{fig:8} and \ref{fig:10}, we present
graphs of the function $S$ which describe the position of the phase change front. The black
and red line represents the analytical and numerical solution respectively.

\newpage${}$
\newpage${}$
\newpage${}$
\newpage
\section{Conclusions}

In this paper, we have developed a numerical method for solving the two-phase fractional
Stefan problem, which can be used as an effective tool for determining a solution alternative to the
closed analytical one. The applied special case of front-fixing technique allows us use the mesh with
constant spatial and time step, which provides discretization of the Caputo derivative with
respect to time. Let us point out that the proposed method has some limitations directly derived
from the construction of the new spatial variables requiring the analytical explicit form of the
function describing the moving boundary depending on the unknown parameter $p$. 

Based on the numerical results presented in Section 5, we can draw the following conclusions:
\begin{enumerate}
	\item The error generated by the numerical solution of the governing subdiffusion equations is 
	a decreasing function of the order of the Caputo derivative. In the special case for $\alpha=1$ 
	,we observe very good matching of the numerical solution to the analytical one.
	\item The relationship between the accuracy of the numerical solution of function $S$ and the order
	of the Caputo derivative is not as obvious as it was listed in the previous point. This is 
	confirmed by the graphs in the Figures \ref{fig:6} and \ref{fig:8}. It should be emphasized 
	that the accuracy of the numerical solution of the function $S$ is the result of many factors, 
	which explains formula (\ref{eq:Sfin}).
	First, the numerical solutions of the functions $u_1$ and $u_2$ are burdened by different
	errors. Second, the approximation of the heat flux calculated in the liquid phase is usually 
	more accurate than that calculated for the solid phase, because for the example considered in 
	the paper, there is the following relationship $x_{m_1,j}-x_{m_1-1,j}< x_{1,j}-x_{0,j}$ (see 
	formula (\ref{eq:Sfin}) and Figure \ref{fig:mesh}). In summary, the accuracy of the numerical solution of 
	the function S is not only dependent on the order of Caputo derivative $\alpha$, but also depends on mesh 
	parameters and the model's physical parameters. 
	\item The proposed numerical method is particularly useful for large values of $\alpha$, i.e. 
	when the closed analytical solution is difficult to calculate (for large values of the similarity 
	variable $\frac{x}{\tau^{\alpha/2}}$).
	
\end{enumerate}

\section*{Acknowledgements} 
This research was supported by the Czestochowa University of Technology Grant Number BS/MN-1-105-301/17/P.
I would like to express my warm thanks to Vaughan R. Voller, who gave scientific guidance that greatly 
assisted the research.

\bibliography{References}%
\bibliographystyle{unsrt}
\end{document}